\documentclass[11pt]{amsart}
\usepackage{amscd}
\newtheorem{thm}{Theorem}

\newtheorem{lem}{Lemma}
\newtheorem{cor}[thm]{Corollary}

\theoremstyle{remark}

\theoremstyle{definition}
\newtheorem{defn}{Definition}

\newcommand{\C}{\mathbb{C}}

\newcommand{\F}{\mathbb{ F}}

\newcommand{\Z}{\mathbb{ Z}}

\newcommand{\R}{\mathbb{ R}}

\newcommand{\id}{\operatorname{id}}

\newcommand{\ra}{{\rightarrow}}

\newcommand{\Fp}{\mathbb{F}_p}

\title{Equivariant Gysin maps and pulling back fixed points}

\author{Bernhard~Hanke}
\address{B.~Hanke, Universit\"at M\"unchen, Theresienstr. 39, 80333
M\"unchen, Germany}
\email{Bernhard.Hanke@mathematik.uni-muenchen.de}

\author{Volker~Puppe}
\address{V.~Puppe, Universit\"at Konstanz, 78457 Konstanz, Germany}
\email{Volker.Puppe@uni-konstanz.de}

\begin{document}

\begin{abstract} We develop a new approach to the 
pulling back fixed point theorem of W.~Browder and use 
it in order to prove various generalizations 
of this result.  
\end{abstract}

\date{\today; MSC 2000: primary 55N20, 55N91, 57S17; secondary 57N65}
\thanks{The first  author is a member of 
the {\sl European Differential Geometry Endeavour} (EDGE), 
Research Training Network HPRN-CT-2000-00101, supported by 
The European Human Potential Programme.}

\maketitle

\section{Introduction}

In \cite{Br}, W.~Browder proves 

\begin{thm} \label{browd} Let $G$ be a finite abelian 
$p$-group, let $M$ be an oriented smooth $G$-manifold and let $N$ 
be an oriented PL $G$-manifold. Assume that $M$ and $N$ are without boundary and 
have the same dimension. Let   
\[
    f : M \ra N
\]
be a continuous proper $G$-equivariant map. Additionally, assume 
that if $p=2$ and $H < G$, then the normal bundle of $M^H$ in $M$ 
has the structure of a complex linear $G$-bundle. Then, if the degree of $f$ 
is not divisible by $p$, the induced map of fixed point sets
\[
    f^G : M^G \ra N^G
\]
is surjective (put differently, any  point in $N^G$ 
can be pulled back under $f$ to a point in $M^G$). 
\end{thm}

Actually, the main result, Theorem (1.1) in \cite{Br}, says that under
the above conditions, the induced map in $\F_p$-cohomology 
\[
  (f^G)^* : H^*(N^G ; \Fp) \ra H^*(M^G ; \Fp)
\]
is injective. The conclusion stated above follows from this fact. 

As a corollary of Theorem \ref{browd}, Browder shows that 
a finite abelian $p$-group cannot act  smoothly on an oriented closed 
manifold with exactly one fixed point (with the extra assumption, 
if $p=2$). This was proved independently by Ewing and Stong 
in \cite{ES} who also characterize those finite groups that can act 
smoothly on a closed manifold with exactly one fixed point. In particular, 
they show that for any finite group which is 
not abelian or contains elements
of coprime order such actions exist (even if one requires 
the underlying manifold to be oriented and the 
action to be orientation preserving).   
A concrete such example with  a group of order $p^3$ 
acting on a complex projective space is also described in Browder's 
paper. In particular,  Theorem \ref{browd} cannot be
extended to a larger class of finite groups.

Theorem \ref{browd} was shown by an involved argument 
relying heavily on the differentiability
of the action on $M$. On the other hand, if $G$ is 
an elementary abelian $p$-group, it is well known 
that this theorem holds more generally for actions on $\Fp$-homology manifolds. 
The proof of this last result uses 
equivariant (ordinary) cohomology and the localization  theorem \cite{AP, Bre}, standard techniques in Smith theory. 

As Browder points out at the end of his paper, it 
would be desirable to prove versions of his
theorem in more general contexts, above all to weaken the differentiability assumption on $M$. 

In this paper, we will provide such  generalizations and 
in particular, if $G$ is a finite cyclic $p$-group, we remove the 
assumption of differentiability of the $G$-action 
on $M$. Our 
approach is different from Browder's. It is based on a 
combination of the Atiyah-Segal-tom Dieck localization 
theorem with equivariant Gysin maps for 
generalized equivariant cohomology theories. 
In this respect, our discussion provides
a link of Browder's result to methods building on classical Smith theory. 

The authors are grateful to W.~Browder for useful 
discussions about the subject treated in this paper and
to  S.~Marde\v si\'c and Y.~Rudyak for helpful comments concerning 
\v{C}ech-extensions of and Poincar\'e duality with 
respect to  generalized (co-)homology theories.  

\section{Main results}

Let $G$ be a finite group and let $E$ be a pointed complex oriented ring spectrum. Then $E$ 
induces homology and cohomology theories defined on $CW$-pairs 
as usual. The cohomology theory $E^*$ has a \v{C}ech extension 
(denoted by $E^*$ again) defined on pairs of arbitrary topological 
spaces (cf.~\cite{MS}) such that the following tautness
assumption is satisfied: If $X$ is a paracompact 
space and $A \subset X$ is a closed subspace, we have
\[
   E^*(X,A) = \lim_{A \subset U} E^*(X,U) \, ,
\]
where $U$ runs over all open neighbourhoods of $A$ in $X$. 

Let $EG \to BG$ be the universal $G$-fibration, where $EG$ is 
considered as a right $G$-space. The Borel construction 
\[
     X \mapsto X_G = EG \times_G X\, , \, \, (f:X \ra Y) \mapsto f_G= 
\id \times_G f : X_G
\ra Y_G \, ,
\]
defines a functor from the category of $G$-spaces and $G$-maps 
to the category of spaces and maps, and so we get  an associated stable 
equivariant cohomology  theory $E_G^*$ (graded 
over the integers) in the sense of \cite{tD1}, Def. (II.6.8), by setting 
\[
 E_G^*(X) = E^*(X_G) \, , \, \, etc. \, . 
\]
Using the complex orientation of $E$, we consider the multiplicative subset  
\[
   S \subset E^*(BG)
\]
consisting of Euler classes of complex $G$-representations without 
trivial direct summand (each complex $G$-representation $V$ gives
rise to a complex vector bundle $EG \times_G V \ra BG$). 

The following version of the localization theorem (which might be 
of independent interest) will be shown by a slight 
modification of the discussion in \cite{tD1}, (III.3). We remark that our proof 
is independent of the equivariant tautness assumption 
stated in \cite{tD1}, Def. (III.3.5), which 
is not clear for our theory $E_G^*$.

\begin{thm} \label{loca}
Let $X$ be a paracompact $G$-space of finite covering dimension. Then 
the inclusion $X^G \hookrightarrow X$ induces an isomorphism
\[
    S^{-1} E^*(X) \cong S^{-1} E^*(X^G) \, .
\]
\end{thm}

Throughout the paper (and if not indicated otherwise), the term $G$-manifold
stands for a closed topological manifold $X$ equipped with a  
topological $G$-action such that $X$ is a $G$-ANR. By 
the well known results from \cite{Ja}, this additional 
requirement is certainly fulfilled, if either the 
action on $X$ is locally linear or $X$ is triangulable and $G$ acts by 
simplicial maps. We remind the reader of the fact that 
the fixed point set of a simplicial group action on a 
topological manifold need not be a topological manifold.  

Let $M$ be a $G$-manifold. By applying the Borel construction 
to the total space $TM = M \times M$ (equipped with the diagonal $G$-action) 
and the base space $M$ of the tangent microbundle  of $M$, we get 
an induced
microbundle $(TM)_G$ over $M_G$ (note that the induced diagonal $G$-actions on $EG \times (M \times M)$, 
resp.~on $EG \times M$ are free).

\begin{defn} \label{orient} By an {\em $E_G$-orientation} of $M$, we 
mean  an $E$-orientation, i.e.~a  (generalized) Thom class, 
of the microbundle 
$(TM)_G$ over $M_G$ in the nonequivariant sense (cf.~\cite{Sw}, Def. 14.5).
\end{defn}

By restriction to the fibre in the Borel construction, it is clear 
that any $E_G$-orientation of $M$ induces an $E$-orientation of $M$.

Let $M$ and $N$ be $E_G$-oriented $G$-manifolds of the same
dimension. Let $f: M \ra N$ be a continuous $G$-map. Then $f$ induces
a (degree preserving) Gysin map $f_!: E^*(M) \ra E^*(N)$ (cf.~Section \ref{Gysin}). 
The element $f_!(1) \in E^0(N)$ (where $1 \in E^0(M)$ is the canonical unit)
is called the {\em $E$-degree} of $f$. Using
Theorem \ref{loca} and the equivariant Gysin map constructed in Section \ref{Gysin}, 
we get the following analogue of Browder's
injectivity result stated at the beginning of this paper.

\begin{thm} \label{main1}   If the $E$-degree $f_!(1)$ is a unit
in $E^0(N)$, then $f$ induces a split injection of 
$S^{-1} E^*(BG)$-modules
\[ 
    S^{-1} E^*(BG \times N^G) \ra S^{-1} E^*(BG \times M^G) \, .
\]
\end{thm}

This theorem contains useful information, only if the 
left hand side is different from zero. This can 
only happen, if $E$ has the property described 
in the next definition. 

\begin{defn} The theory $E$ is called {\em $G$-sensitive}, if
\[
   S^{-1} E^*(BG) \neq 0 \, .
\]
\end{defn}

Typical examples of $G$-sensitive theories are  $\check{\rm C}$ech cohomology 
with $\Fp$-coefficients
for elementary abelian $p$-groups \cite{AP}, 
complex $K$-theory for cyclic $p$-groups (Lemma \ref{sensitive}) and unitary 
bordism for finite abelian $p$-groups \cite{tD2}. On the 
other hand $\check{\rm C}$ech cohomology  with integer or 
$\Fp$-coefficients is not $\Z / p^2$-sensitive  and complex $K$-theory is not $G$-sensitive, 
if  $G = \Z/p \times \Z/p$ or  $G = \Z / n$, if $n>1$ is 
not a prime power. For non abelian $G$, none of these theories is
$G$-sensitive (cf.~\cite{tD2}, proof of Satz 3).

\begin{thm} \label{main2}  Let $M$ and $N$ be $E_G$-orientable $G$-manifolds
of the same dimension, where $E$ is $G$-sensitive. Furthermore, let
$M$ and $N$ be oriented  in the usual sense and let 
$f : M \ra N$ be a continuous $G$-map of ordinary 
degree $d \in \Z$. If $d  \cdot 1 \in E^0$ is a unit,  then the induced map of fixed points
\[
     f^G : M^G \ra N^G
\]
is surjective. 
\end{thm}

Theorem \ref{main2} gains its power by judicuous choices of 
the theory $E$. Using unitary bordism, we essentially reprove
Browder's result. On the other hand, with the use of 
complex $K$-theory, we 
obtain the following generalization of Browder's theorem 
in the case of actions of cyclic $p$-groups.

\begin{thm} \label{general2} Let $G$ be a finite cyclic $p$-group, where 
$p$ is odd. Let $M$ and $N$ be  oriented (in the usual sense) $G$-manifolds of the same dimension. If $f:M \ra N$ is a $G$-map whose 
(nonequivariant) degree is not a multiple of $p$, then the induced map of fixed point sets $f^G: M^G \ra N^G$ is surjective. 
\end{thm} 

We remark that this theorem  would hold for any finite abelian $p$-group ($p$ odd), if the 
following conjecture is true.

\medskip

\noindent {\bf Conjecture.}  There exists a complex oriented
cohomology theory $E^*$ with the following properties: 
Every closed oriented manifold is orientable with respect to $E^*$
and $E^*$ is $G$-sensitive for every finite abelian $p$-group $G$, where $p$ is an odd prime.
\medskip

If $E$ were a cohomology theory with these properties, it would follow from the proof of 
Lemma \ref{Sullivan} that if $X$ 
were a $G$-manifold which was orientable in the usual sense, then 
$X$ would also be $E_G$-orientable.
A natural canditate for a theory described in this conjecture 
is topological bordism, for which the first property is immediate.
Unfortunately, it seems to be beyond present knowledge, if 
this theory is also $G$-sensitive for  every finite abelian $p$-group $G$ (where
$p$ is odd). 
Nevertheless, it is obvious that topological bordism is $G$-sensitive, 
if $G$ is a cyclic $p$-group 
or if $G$ is an elementary abelian $p$-group:  We
have both Sullivan's $K$-theory orientation and the ordinary 
homology orientation of the oriented bordism spectrum localized
at an odd prime. This gives some indication that the above conjecture
might indeed be true. 

The difficulty that a theory is not $G$-sensitive can sometimes
be circumvented by an inductive procedure based on the ideas of the 
proof of Theorem \ref{main1}. This will be illustrated  
in the last part of this paper, where we sketch a third approach 
to Browder's result using classical Smith 
theory. Moreover, assuming some strong cohomological 
restrictions on the spaces involved, we even can 
remove the assumption that $G$ is abelian: 

\begin{thm}  Let $M$ and $N$ be compact 
oriented $\Fp$-cohomology manifolds of the same dimension (we use \v{C}ech-cohomology), let $G$ be 
a finite $p$-group and assume the following: The 
$\Fp$-cohomology of $M$ and $N$ is concentrated in even 
degrees and 
\[
    p > {\rm max} \{ \dim_{\Fp} H^*(M;\Fp) , \dim_{\Fp} H^*(N;\Fp) \} \, .
\]
Let $f:M \ra N$ be $G$-equivariant. If the degree 
of $f$ is a nonzero number in $\Fp$, then the induced
map $M^G \ra N^G$ is surjective. 
\end{thm}

\section{Localization theorem} 

In this short section, we will prove Theorem \ref{loca} above. 
Our proof will be modelled on the proof of Theorem (III.3.6) in 
\cite{tD1}. We set $A=X^G$. At first, notice that by our assumptions and because 
$G$ is finite, $X - A$ is of finite $S$-type in the sense of \cite{tD1}, 
Def. (III.3.2). We now have to show that 
\[
   S^{-1} E_G^*(X,A) = 0 \, .
\]
Using the (nonequivariant) tautness of $E^*$, it is enough 
to show that 
\[
   S^{-1} E^*(X_G, U) = 0 \, , 
\]
for all open $A_G \subset U \subset X_G$. At this point, we remark
that the projection $X_G \ra BG$ restricts to a map $U \ra BG$ so 
that $E^*(U)$ has indeed the structure of an $E^*(BG)$-module. By excision
\[
   S^{-1} E^*(X_G, U) \cong S^{-1} E^*(X_G - A_G , U - A_G) \, .
\]
Because $X - A$ is of finite $S$-type, we have
\[
   S^{-1} E^*((X-A)_G) = 0\, 
\]
by \cite{tD1}, Theorem (III.3.3), and because the inclusion $U - A_G \hookrightarrow (X-A)_G$ 
induces a map of $E^*(BG)$-algebras $E^*((X-A)_G) \ra E^*(U - A_G)$, 
this implies that also 
\[
   S^{-1} E^*(U - A_G) = 0 \, .  
\]
The long exact cohomology sequence gives 
\[
  S^{-1} E^*(X_G - A_G , U - A_G) = 0  \, ,  
\]
and this completes the proof of Theorem \ref{loca}. 

A relative version of Theorem \ref{loca} follows
easily by use of the long exact cohomology 
sequence. It is also obvious that we could have
proven a version of Theorem \ref{loca} for 
infinite groups $G$, but we restricted
ourselves to the case treated here in view
of our later applications.  

\section{Equivariant Gysin maps} \label{Gysin}

Let $G$ be a finite group, let $M$ and  $N$ be $E_G$-oriented $G$-manifolds
of dimension $m$ and $n$ respectively. Let $f : M \ra N$ be a continuous 
$G$-equivariant map. We will construct an equivariant 
Gysin map 
\[
   (f_G)_!:  E_G^{k}(M) \ra E_G^{k+ n -m}(N)\, , \, k \in \Z\, ,
\]
along the lines of \cite{Ka}. By \cite{Mo},  there is a $G$-embedding 
$M \ra V$ into a complex $G$-representation. Let $c$ be its complex 
dimension. Using the map $f$ and this embedding, we consider $M$ as a $G$-submanifold 
of $N \times V$ (equipped with the diagonal $G$-action). 
In a first step, we need to construct a Thom-like class
\[
  \tau \in E^{n + 2c - m}((N \times V)_G, (N \times V)_G - M_G) \,
\]
which is straightforward, if $M$ has a (stable) normal $G$-microbundle 
in $N \times V$. Because the literature (cf.~\cite{Ib}) provides
such a bundle, only if the $G$-spaces under consideration 
are locally linear, we proceed as follows. 

As $M$ is a $G$-ANR (by our definition of $G$-manifold), there is a 
$G$-invariant neighbourhood
$U \subset N \times V$ of $M$ and a $G$-equivariant retraction
\[
    r: U \ra M \, .
\]
We choose a $G$-CW model of $EG$  and obtain a skeletal filtration
\[
   EG^0 \subset EG^1 \subset EG^2 \subset \ldots
\]
of $EG$. Write $M_G^i$ for $EG^i \times_G M$ and so forth. At first, we  
construct a compatible sequence of Thom classes 
\[
    \tau^i \in E^{n + 2c - m}((N \times V)_G^i, (N \times V)_G^i - M_G^i)\, .
\]
Using the $E_G$-orientation of $N$ and the canonical $E_G$-orientation of $V$, we get an $E_G$-orientation of the $G$-manifold $N \times V$. By 
restriction, we then get $E$-orientations of the microbundles
$(TM)^i_G$ over $M_G^i$ and of $(T(N \times V))_G^i$ over
$(N \times V)_G^i$. 

We will construct a normal microbundle $\nu$ of $M_G^i$ in $(N \times V)_G^i \times \R^q$ for some $q$  along the lines of \cite{Mi}. 
By this, we mean that (the total space of) $\nu$ is homeomorphic
to a neighbourhood of $M_G^i$ in $(N \times V)_G^i \times \R^q$ 
by a homeomorphism mapping the zero section of $\nu$ to $M_G^i$. 

The retraction $r$ induces a retraction
\[
    r_G^i : U_G^i \ra M_G^i \, .
\]
Denoting a microbundle and 
its total space with the same letter, the proof of \cite{Mi}, Lemma (5.3), 
can be carried out fibrewise in the Borel fibration and shows 
that we have a homeomorphism 
\[
    (TU)_G^i~|_{M_G^i} \approx (r_G^i)^*(TM)_G^i \, . 
\]
In particular, the total space of $(TU)_G^i|_{M_G^i}$ constitutes a microbundle 
neighbourhood (defined in the obvious sense) of $M_G^i$  in $(r_G^i)^*(TM)_G^i$ 
(cf.~\cite{Mi}, Lemma (5.4)). 
The space $M_{G}^i$ being compact, there is a stably unique microbundle
$\eta$ over $M_G^i$ complementary to $(TM)_G^i$, i.e.~its Whitney sum with $(TM)_G^i$ is a trivial 
bundle $\underline{\R}^q$. As in 
\cite{Mi}, Theorem (5.8), the bundle $\eta$ can be used in order to construct a 
normal microbundle $\nu$ of the required sort. 

Theorem (5.9) in \cite{Mi} can also be adapted to our situation 
and shows that we get a microbundle isomorphism 
\[
   (TM)_G^i \oplus \nu \cong (T(N \times V))_G^i~|_{M_G^i} \oplus 
\underline{\R^q} \, .
\]
For a translation of  the last part of the proof of Theorem (5.9)
in \cite{Mi}, choose a $G$-invariant neighbourhood $D$ of the diagonal in
 $M\times M$ such that $p_1|_{D}$ and $p_2|_{D}$ (where $p_1$ and 
$p_2$ are projections) are equivariantly
homotopic (this can be done because $M$ is a $G$-ANR) and 
apply the Borel construction to $D$. 

By \cite{Dyer}, the microbundle $\nu$ has a canonical $E$-orientation. 
The class $\tau^i$ is the $q$-fold desuspension of the Thom class of this bundle. 

We want to show that under the inclusion $M_G^i \subset M_G^{i+1}$, 
the Thom class $\tau^{i+1}$ is mapped to $\tau^i$. But this is obvious, 
because a complementary bundle to  $(TM)_G^{i+1}$ restricted
to $M_G^i$ is complementary to $(TM)_G^i$ and the rest of
the construction of the normal bundle $\nu$ above was based
on data which were independent of $i$. 

Because the canonical map 
\begin{eqnarray*}
  \lefteqn{\rho:   E^{n + 2c - m}((N \times V)_G, (N \times V)_G - M_G)
        \to } \hspace{4cm} \\
     & & \lim_i E^{n + 2c - m}((N \times V)_G^i, (N \times V)_G^i - M_G^i)  
\end{eqnarray*}
in the Milnor exact sequence is surjective, we can choose a  class 
\[
  \tau \in E^{n + 2c - m}((N \times V)_G, (N \times V)_G - M_G)
\]
restricting to the inverse system $(\tau^i)$. But notice that due to the possible nonvanishing of $\lim^1$, 
this class need not be uniquely determined by $(\tau^i)$. 

Multiplication with $\tau$ from the right induces a map
\begin{eqnarray*}
  \lefteqn{\phi:  E^k(M_G) \stackrel{r_G^*}{\ra} E^k(U_G) \stackrel{\cdot \tau}{\ra} E^{k+n+2c -m}(U_G, U_G - M_G) 
\cong} \hspace{4cm} \\
     & &  E^{k+n+2c -m}((N \times V)_G, (N \times V)_G - M_G) \, .
\end{eqnarray*}
This is an isomorphism, because it induces isomorphisms 
\[
    E^k(M_G^i) \cong E^{k+n+2c -m}((N \times V)_G^i, ((N \times V)_G^i - M_G^i) 
\]
for all $i$ (and hence induces isomorphisms of the $\lim^0$ and 
$\lim^1$ terms in the respective Milnor exact sequences). In 
this respect, we might consider $\tau$ as a Thom class.

Now choose a $G$-invariant disc $D(V) \subset V$ such that $M \subset N \times D(V)$ and denote  its boundary by $S(V)$.
The equivariant Gysin map $(f_G)_!$ is defined as the composition
\begin{eqnarray*}
   E^k(M_G) & \stackrel{\phi}{\ra} & E^{k+n+2c-m}((N \times V)_G, (N \times V)_G - M_G) \\
            & \ra & E^{k+n+2c-m}((N \times D(V))_G , (N \times S(V))_G) \\
       & \ra & E^{k+n - m}(N_G)
\end{eqnarray*}
where the first map is the Thom isomorphism constructed before, the second map is induced by an inclusion of spaces and the last map is 
the inverse of the Thom isomorphism of the complex 
vector bundle
\[
     V \hookrightarrow (N \times V)_G \ra N_G \, .
\]

The following properties of $(f_G)_!$ are proven similarly 
as in \cite{Ka}. 
\begin{itemize}
   \item Under restriction to the fibres in the Borel constructions,
         $(f_G)_!$ restricts to the usual Gysin map 
    \[
        f_! : E^k(M) \ra E^{k+n-m}(N)
    \]
        where we use the induced $E$-orientations of $M$ and $N$. 
        Recall that $f_!$ is induced by the map in homology $f_*: E_*(M) 
        \ra E_*(N)$  and the Poincar\'e duality isomorphisms for $M$ and $N$.
   \item $(f_G)_!$ is an $E^*(BG)$-module map.
   \item For $x \in E_G^*(N)$, we have 
        \[
           (f_G)_!(f_G^*(x)) = x \cup (f_G)_!(1) \, .
        \]
\end{itemize}

From now on, we assume that $m = n$.

\begin{lem} \label{iso} If the $E$-degree $f_!(1) \in E^0(N)$ is a unit,
then the composition
\[
   E_G^*(N) \stackrel{f_G^*}{\ra} E_G^*(M) \stackrel{(f_G)_!}{\ra} E_G^*(N)
\]
is an isomorphism. 
\end{lem}

\begin{proof} We use a 
spectral sequence argument. The Borel fibrations 
for $M$ and $N$ lead to converging Atiyah-Hirzebruch spectral sequences
(\cite{Sw}, Theorem 15.7.) 
\[
   E_2^{p,q}(M) = H^p(BG ; E^q(M)) \Longrightarrow E_G^{p+q}(M)
\]
and similarly for $N$. Clearly, the map $(f_G)^*$ induces
a map of $E_2$-terms. But also $(f_G)_!$ induces
a map of $E_2$-terms
\[
  H^p(BG ; E^q(M)) \rightarrow H^p(BG ; E^q(N)) 
\]
which is induced by the non equivariant Gysin map $f_!:E^q(M)  \ra E^q(N)$.
Hence, on the $E_2$-level, the composition $(f_G)_! \circ (f_G)^*$ 
is induced by the map of coefficients
\[
   E^q(N) \stackrel{f^*}{\ra} E^q(M) \stackrel{f_!}{\ra} E^q(N)\, .
\]
This map is given by $x \mapsto x \cup f_!(1)$ and hence
is an isomorphism by our assumption. Consequently, 
the induced map of $E_{\infty}$-terms is an isomorphism, too. This implies
our assertion.  
\end{proof}

In order to  start with the proof of Theorem \ref{main1},  we consider
the multiplicative subset $S \subset E_G^* = E^*(BG)$ generated
by Euler classes of irreducible complex $G$-representations without
trivial direct summand. Theorem \ref{loca} yields a commutative localization
diagram 
\[
  \begin{CD}
      E_G^*(N)    @>f_G^*>>    E_G^*(M)     @>(f_G)_!>> E_G^*(N) \\
     @VVV                     @VVV          \\
S^{-1} E_G^*(N) @>>>      S^{-1}E_G^*(M)  \\
     @V\cong VV               @V\cong VV             \\   
S^{-1} E^*(BG \times N^G)  @>({\rm id} \times f^G)^*>>     S^{-1}E^*(BG \times M^G) 
  \end{CD}
\] 
The first line is an isomorphism by Lemma \ref{iso} and 
in particular, the map $f_G^*$ is a split injection of 
$E^*(BG)$-modules. This 
is therefore also true for the map 
$(\id \times f^G)^*$ in the third line. Hence, the proof
of Theorem \ref{main1} is complete.

\medskip

Now additionally assume that $M$ and $N$  are oriented in 
the usual sense.  
  
\begin{lem} \label{unit}  Let $d \in \Z$ be the usual degree 
of $f$. If $d \cdot 1 \in E^0$ is a unit, then also $f_{!}(1) \in E^0(N)$ is a unit.
\end{lem}

\begin{proof} If $X$ is a closed oriented topological manifold
of dimension $n$, obstruction theory yields a 
canonical bijection between the set $[X,S^n]$ of homotopy classes of maps $X \ra S^n$ and elements of $H^n(X;\Z)$ 
(recall that $X$ is homotopy equivalent 
to a $CW$-complex). Let $\pi:M \ra S^n$ and $\psi:N \ra S^n$ be 
maps of degree $1$ (given by collapsing the exterior of a small ball to a point). 
If we take a map $S^n \ra S^n$ of degree $d$, we hence get a diagram  
\[
  \begin{CD}   
            M         @>f>>       N         \\
        @V\pi VV               @V\psi VV      \\
           S^n    @>h>>          S^n         \\
  \end{CD} 
\]   
commuting up to homotopy. We can choose orientation classes $[M] \in \tilde{E}_n(M)$ and 
$[N] \in \tilde{E}_n(N)$ such that $\pi_*([M]) = [S^n] = \psi_*([N])$, where 
$[S^n] \in \tilde{E}_n(S^n)$ denotes the canonical orientation class (i.e.~the $n$-fold suspension of
$1 \in \tilde{E}_0(S^0)$). In particular, we have $(h \circ \pi)_*[M] = d \cdot [S^n]$.  
Because $(h \circ \pi)_!(1)$ is the Poincar\'e dual of $(h \circ \pi)_*([M])$, 
we get the equation
\[
     \psi_! \circ  f_! ( 1) = d \cdot 1 \in E^0 \subset E^0(S^n) \, .
\]
We have canonical splittings 
\[
   E^0(N) = \tilde{E}^0(N) \oplus E^0\, , \, \, E^0(S^n) = \tilde{E}^0(S^n)
                           \oplus E^0 \, 
\]
and by Alexander duality  (\cite{Sw}, Theorem 14.11) we see that $\psi_!$ induces a map
of reduced groups $\tilde{E}^0(N) \ra \tilde{E}^0(S^n)$. Altogether,
this implies that the component of $f_!(1)$ in $E^0$ must be 
equal to $d \cdot 1$. 
Because $N$ is finite dimensional, $\tilde{E}^0(N)$ is nilpotent (apply
the Atiyah-Hirzebruch spectral sequence) and the element $f_!(1)
\in E^0(N)$ is indeed a unit. 
\end{proof}

For the proof of Theorem \ref{main2}, choose a point $y \in N^G$ and  assume that 
$f^{-1}(y) \cap M^G = \emptyset$.  This implies 
that the composition 
\[
 \chi : S^{-1} E^* (BG \times (N^G , N^G - \{y \} )) \ra S^{-1} E^*(BG \times N^G ) \ra S^{-1} E^*(BG \times M^G ) 
\] 
is zero. 

Let $\pi:N \ra \{y\}$ be the unique $G$-map and consider the composition of $E^*(BG)$-linear maps
\[
 E^*(N_G, (N-\{y\})_G) =   E_G^*(N, N - \{y\}) \ra E_G^*(N) \stackrel{(\pi_G)_!}{\to} E_G^*(\{ y \} ) = E^*(BG) \, ,
\]
where $(\pi_G)_!$ is the equivariant Gysin map considered before. This composition 
decreases the degree by $n$. The pair $(N_G, (N - \{y \})_G)$ can be considered as the Thom space of the restricted bundle  $(TN|_{\{y\}})_G$. Because $N$ is $E_G$-orientable, this 
restricted bundle is $E$-orientable and the composition 
under consideration is an isomorphism: It follows directly from the construction of Gysin maps at the beginning of section \ref{Gysin} that for each $i$, the corresponding map 
\[
    E^*(N_G^i, (N-\{y\})_G^i) \ra E^*(BG^i)
\]
is exactly the inverse of the Thom isomorphism 
of the $E$-oriented bundle $(TN|_{\{y\}})_G^i$ and 
therefore is an isomorphism.  An application 
of the Milnor exact seqence to the induced isomorphism 
of inverse systems shows that the composition 
\[
  E^*(N_G, (N-\{y\})_G) \ra  E^*(BG)
\]
is an isomorphism, as well. In particular, the induced
localized map  
\[
   S^{-1} E^*(N_G, (N-\{y\})_G) \ra S^{-1} E^*(BG)
\]
is surjective. Its image is different from zero, because $E$ is $G$-sensitive. Application of the localization theorem  finally 
shows that the image of the first map in  the composition 
$\chi$ is not zero. 

By Lemma \ref{unit} and 
Theorem \ref{main1},  the second map in this composition is injective. 

Hence, we get a  contradiction to  the fact that the composition 
$\chi$ is zero.  The 
assumption $f^{-1}(y) \cap M^G = \emptyset$ must therefore be 
false and the proof of Theorem  \ref{main2} is complete.

\section{Applications and examples}

The first choice for $E$ that comes to mind is \v{C}ech cohomology with 
$\Fp$-coefficients.  If $G$ is an elementary abelian $p$-group,
$M$ and $N$ are closed $G$-manifolds which are oriented, 
if $p$ is odd, and 
if $f : M \ra N$ is a $G$-equivariant map  
whose degree (as an element in $\Fp$) is different from $0$, then all the 
assumptions stated in Theorem \ref{main2} are fulfilled and the induced map of fixed point 
sets is surjective.  It is well known that this  conclusion
holds as well if $M$ and $N$ are just assumed to 
be compact oriented $\Fp$-homology manifolds.  Also, note that the conclusion of Theorem \ref{main1} holds, if we just assume that $M$ and $N$ are compact $\Fp$-Poincar\'e 
duality spaces. In particular, the induced map $H^*(N^G ; \Fp) \ra H^*(M^G; \Fp)$ 
is injective in this case.

Now let $G$ be a finite abelian $p$-group. For $E$, we take 
the unitary bordism spectrum localized at $p$ (and denoted by $MU_{(p)}$). It
follows from results of tom Dieck \cite{tD2} that $MU_{(p)}$ is 
$G$-sensitive. 

\begin{lem} \label{MU} Let $p$ be odd and let $M$ be an oriented closed smooth 
$G$-manifold. Then $M$ is $(MU_{(p)})_G$-orientable.
\end{lem}

\begin{proof} For odd $p$, we have an orientation (i.e. a map of ring spectra)
\[
    {\rm MSO}_{(p)} \ra MU_{(p)} \, , 
\]
c.f.~\cite{MM} or \cite{Ru}.  In particular, any oriented vector bundle  
is $MU_{(p)}$-orientable. If $M$ is an oriented smooth 
$G$-manifold, the induced bundle  $(TM)_G \ra M_G$ is orientable and therefore $MU_{(p)}$-orientable.
\end{proof}

If $p = 2$, one shows similarly that $M$ is $(MU_{(p)})_G$-orientable, if
$M$ is a stably almost complex manifold and the group action 
preserves the stable almost complex structure of $M$. 
Using Lemma \ref{MU} and Theorem \ref{main2} applied to 
$E= MU_{(p)}$, we get the following version of Browder's original theorem.

\begin{thm} Let $G$ be a finite abelian $p$-group. Let $M$ 
and $N$ be closed oriented smooth $G$-manifolds 
if $p$ is odd, or closed stably complex $G$-manifolds if $p=2$. 
Assume that $M$ and $N$ have the same dimension. 
If $f : M \ra N$ is an equivariant map whose degree is not 
divisible by $p$, then the induced map $M^G \ra N^G$ is 
surjective. 
\end{thm} 

Note that for reproving 
Browder's theorem (in the differentiable case), we need not 
refer to the construction of equivariant Gysin maps in Section \ref{Gysin},
but can use directly the corresponding results from \cite{Ka}.

The following well known  example (cf.~\cite{ES}) shows that the case 
$p=2$ really requires special attention: The action of $\Z/4$ on $\C P^2$ 
induced by 
\[
    [x_0:x_1:x_2] \mapsto [\overline{x_0}: - \overline{x_2}: 
\overline{x_1}]
\]
is orientation preserving, smooth and has exactly one fixed point.

Now let $G$ be a finite cyclic $p$-group. As our theory, 
we consider $K_{(p)}$, complex $K$-theory localized at $p$. 

\begin{lem} \label{sensitive} Komplex $K$-theory localized at $p$ is $G$ sensitive. 
\end{lem} 

\begin{proof} (cf.~\cite{tD1}, p.~33 f.) Let $G = \Z / p^r$, $r \geq 1$. By application 
of the Gysin sequence, we have an isomorphism
\[
   K^*_{(p)}(BG) \cong K^*_{(p)}[[C]]/(e(\eta^{p^r}))
\]
where $C \in K^2_{(p)}(\C P^{\infty})$ is the Euler class of the universal 
line bundle on $\C P^{\infty}$ and $e(\eta^{p^r}) \in K^2_{(p)}[[C]]$ is the Euler class of the $p^r$th 
tensor power of this line bundle.  One sees that $e(\eta^{p^r})$ has 
a factorization containing $C$ and each of the power series
\[
    e(\eta^{p^i})/ e(\eta^{p^{i-1}})\, , \,  1 \leq i \leq r \, ,
\]
With help of the formal group law for complex $K$-theory
\[
    (x,y) \mapsto x + y + xy \, ,
\]
these power series factors are all prime (because they all have constant term $p$) 
and different (because they are polynomials
in $C$ of different degrees). In particular, $e(\eta^{p^r})$ 
contains a prime factor not contained in any $e(\eta^{p^i})$, $0 \leq i < r$.
(Here we use the fact that $K^*_{(p)}[[C]]$ is a unique factorization 
domain.)  
The multiplicative set $S$ is generated by the Euler classes
\[
   e(\eta^{p^i}) \, , \, 0 \leq i < r\, ,
\]
hence, by the argument above, $S \cap (e(\eta^{p^r})) = \emptyset$. 
This implies that the localization $S^{-1} K^*_{(p)}(BG)$ 
is indeed different from zero. 
\end{proof}

By a similar argument, one also sees that $K^*_{(p)}$ is 
not $\Z/n$-sensitive, if $n > 1$ is divisible by two different 
primes $p$ and $q$: In the 
localization $S^{-1} K^*_{(p)}(\Z/n)$, certain power series in $C$ starting
with constant terms $p$ and $q$  are 
zero and therefore this is also true for some power series starting with $1$.
But such a power series is a unit
in $K^*_{(p)}[[C]]$.

We now assume additionally that $p$ is odd. 

\begin{lem} \label{Sullivan} Each  oriented (in the usual sense) $G$-manifold $X$ is $(K_{(p)})_G$-orientable.
\end{lem}

\begin{proof} We prove at first that $(TX)_G$ is ${\rm MSTOP}$-orientable. 
Let $(TX)_G^i$ and $X_G^i$ be filtrations of $(TX)_G$ and 
$X_G$ by finite dimensional subcomplexes as before. By Kister's 
theorem \cite{Ki} we can replace the microbundle $(TX)_G^i$ 
over the Euclidean neighbourhood retract $X_G^i$ by an essentially unique 
fibre bundle and hence, for each $i$, we get a tautological Thom class 
\[
     \tau^i \in [M( (TX)_G^i ), {\rm MSTOP}] \, ,
\]
(which is an ${\rm MSTOP}$-cohomology class of the Thom 
space of  $(TX)_G^i$). The Thom classes constructed in this way
are compatible under restriction and so $(\tau^i)$ defines 
an inverse system, which induces a Thom 
class for $(TX)_G$. Sullivan proved that there is an orientation 
(cf.~\cite{MM} or \cite{Ru})
\[
   {\rm MSTOP}_{(p)} \ra K_{(p)} \, 
\]
of topological bordism localized at $p$. This completes the 
proof of Lemma \ref{Sullivan}.
\end{proof} 

Lemma \ref{sensitive} and \ref{Sullivan} together with Theorem \ref{main2} complete the proof of Theorem \ref{general2}.

With regard to Theorem \ref{main1} and the remark following 
Theorem \ref{browd} stated at the beginning of 
this paper, one might wonder, if under the assumptions of 
Theorem \ref{main1} and in the case that $E$ is $G$-sensitive, one can prove 
that the induced map
\[
    E^*(N^G) \ra E^*(M^G)
\]
is injective. However, even if we 
assume that (for compact spaces $X$) we have a K\"unneth formula
\[
    S^{-1} E^*(BG \times X) \cong E^*(X) \otimes_{E^*} S^{-1}E^*(BG)
\]
such a result is not immediate for general theories $E$, because a non 
injective map can very well become an injection after tensoring 
with some identity map.

\section{Iterative use of localization} \label{inductive}

Classical Smith theory allows the discussion of 
Browder's result in the contex of  Poincar\'e duality spaces.  
We use $\check{\rm C}$ech cohomology with $\Fp$-coefficients. Let $M$ and $N$ be compact $\Fp$-Poincar\'e duality spaces of the same dimension on which a finite 
(not necessarily abelian) $p$-group is acting. Furthermore, let $f:M \ra N$ a $G$-equivariant
 map whose degree (as an element in $\Fp$) is 
different from $0$. The ultimate generalization of Browder's 
theorem would be the conclusion that under these assumptions, 
the induced map 
\[
    H^*(N^G) \ra H^*(M^G)
\]
is injective (which implies that $f^G : M^G \ra N^G$ is surjective, 
if the space $N$ has good local cohomology properties,
e.g.~if $N$ is a $\Fp$-cohomology manifold). 
Of course, this 
is known to be false in general, if we do not assume 
$G$ to be abelian. On the other 
hand, in view of Theorem \ref{main2}, this conclusion is 
indeed valid (and classical), if $G$ is an elementary abelian $p$-group. 
But ordinary cohomology (even with integer coefficients) 
is not $\Z / p^2$-sensitive, so we will use 
induction on the order of $G$ (as 
Browder does in his proof of Theorem \ref{browd})  and wish to formulate conditions 
under which the induced map $ H^*(N^G) \ra H^*(M^G)$ is 
injective.

The essence of the argument can be seen in the special case $G= \Z /p^2$.  Let 
$H < G$ be the subgroup
of index $p$ and set $K = G / H \cong \Z /p$. We will 
use some of the result in  \cite{AP}, Chapter 5.  

Considering $f$ as an $H$-equivariant map, one gets 
an equivariant Gysin homomorphism (which is degree preserving as $m=n$) 
\[
    (f_H)_! : H_H^*(M) \ra H_H^*(N) \, ,
\]
where as usual the subscript indicates taking 
cohomology of the associated Borel constructions. Let 
$S \subset H^*(BH)$ (resp.~$T \subset H^*(BK)$) be the multiplicative subset generated 
by a nonzero element in $H^2(BH)$ (resp.~by an nonzero element in $H^2(BK)$). We then get a commutative 
localization 
diagram 
\begin{tiny}
\[
\begin{CD}
H_H^*(N) @>(f_H)^* >> H_H^*(M) @> (f_H)_! >> H_H^*(N) \\
  @VVV                   @VVV                            @VVV   \\
S^{-1} H_H^*(N) @>S^{-1} (f_H)^* >> S^{-1} H_H(M) @> S^{-1} (f_H)_! >> S^{-1} H_H(N)  \\
  @V \cong VV             @V \cong VV          @V \cong VV     \\
H^*(N^H) \otimes S^{-1}H^*(BH) @>(f^H)^* \otimes {\rm id}>> H^*(M^H) 
\otimes S^{-1} H^*(BH) @>\alpha_H>> H^*(N^H) \otimes S^{-1} H^*(BH) \, , 
\end{CD}
\]
\end{tiny}
where $\alpha_H$ is the composition (${\rm mult}$ denotes
multiplication from the right)
\[
    {\rm mult}_{e_H^N}  \circ ( (f^{H})_! \otimes {\rm id}) \circ {\rm mult}_{(e_H^M)^{-1}}
\]
using the equivariant Euler classes 
\[
   e_H^N \in H^*(N^H \times BH)\, , ~~ e_H^M \in H^*(M^H \times BH)
\]
of the inclusions $M^H \subset M$ and
$N^H \subset N$ (these Euler classes  become invertible after localizing) and the nonequivariant Gysin map for the map $f^H$ between $\Fp$-Poincar\'e 
duality spaces. Our assumption 
implies that the first line of this diagram is given by multiplication 
with a nonzero number 
in $\Fp$ and the same holds for the last line showing that 
$(f^H)^*$ is injective.  

For the next step, note that both $N^H \times BH$ and $M^H \times BH$  carry induced $K$-actions (trivial on the classifying space
factors). This suggests taking the corresponding Borel constructions 
and arguing as before (starting with the last line 
of the above diagram). In order to make this 
work, we will assume for the moment that the equivariant Euler classes 
\[
e_H^M \in  H^*(M^H \times BH) \, , ~~e_H^N \in H^*(N^H \times BH)
\]
lie in the image of the  maps
\begin{eqnarray*}
  H^*((M^H)_K \times BH) & \ra &  H^*(M^H \times BH) \, , \\
    H^*((N^H)_K \times BH) & \ra &  H^*(N^H \times BH)
\end{eqnarray*}
induced by restricting to a fibre in the Borel construction.  
We remark that this lifting assumption is not automatic for the following reason: 
The $H$-equivariant Euler classes $e_H^M$ and $e_H^N$ depend on  
neighbourhoods of $M^H$ in $M$ and $N^H$ in $N$. But a priori, 
an action of the quotient group $K$ can only be defined on 
$M^H$ and $N^H$ and not on a neighbourhood of these subsets. 
However, in some cases this assumption can be verified. We will 
discuss this point further below. 

Let $\tilde{e}_H^N$ and $\tilde{e}_H^M$ be preimages 
$e_H^N$ and $e_H^M$ under the maps above. Using a spectral sequence 
argument, one  shows that these classes become again invertible after
inverting $S$. We define 
a map 
\[
  \alpha_H':  H^*((M^H)_K \otimes S^{-1} H^*(BH) \ra  H^*((N^H)_K) \otimes S^{-1} H^*(BH)
\]
as the composition 
\[
    {\rm mult}_{\tilde{e}_H^N}  \circ  ((f^{H})_! \otimes {\rm id})
 \circ {\rm mult}_{(\tilde{e}_H^M)^{-1}} \, .
\]
Using the abbreviations $H'$ for $S^{-1} H^*(BH)$, and $K'$
for $T^{-1} H^*(BK)$, we then 
get a commutative diagram 
\begin{tiny}
\[
\begin{CD}
H^*((N^H)_K) \otimes H' @>(f^H)_K^* \otimes {\rm id}>> 
H^*((M^H)_K) \otimes H' @>\alpha_H'>> H^*((N^H)_K) \otimes H' \\
  @VVV                   @VVV                            @VVV   \\
T^{-1}H^*((N^H)_K) \otimes H' @>T^{-1} (f^H)_K^* \otimes {\rm id}>> T^{-1}H^*((M^H)_K) \otimes H' @>\alpha_H'>> T^{-1} H^*((N^H)_K \otimes H' \\
  @V \cong VV             @V \cong VV          @V \cong VV     \\
H^*(N^G) \otimes K' \otimes H' @>(f^G)^* \otimes {\rm id}>> H^*(M^G) \otimes K' \otimes H' @>\alpha_G >> 
H^*(N^G) \otimes K' \otimes H' 
\end{CD}
\]
\end{tiny}
with 
\[
  \alpha_G' : {\rm mult}_{\tilde{e}_H^N \cup e_K^N } \circ 
            ( (f^H)_! \otimes {\rm id} )\circ 
            {\rm mult}_{(\tilde{e}_H^M \cup e_K^M)^{-1}}
\]
and equivariant Euler classes $e_K^N$, $e_K^M$ of 
the embeddings of $K$-spaces $N^G \subset N^H$ and $M^G \subset M^H$. 
Using a spectral sequence argument again, the first 
line of the last diagram is given by multiplication 
with a nonzero element in $\Fp$ and the induced
map $H^*(N^G) \ra H^*(M^G)$ is indeed injective. 

Note that the property of $G$ being abelian was
not used in the course of this argument. Everything depends
on the lifting assumption above which is fairly obvious in 
the following special case.

\begin{lem} \label{strong} 
Assume that the induced actions of $K$ on $M^H$ and 
on $N^H$ are totally nonhomologous to zero. Then 
the above lifting assumption for the relevant 
$H$-equivariant Euler 
classes is fullfilled (because the relevant 
maps are surjective  in this case). 
\end{lem}

By imposing strong cohomological retrictions on 
the actions involved and  using the fact 
that every finite $p$-group is solvable, one can 
therefore prove

\begin{thm} \label{coho} Let $M$ and $N$ be compact 
oriented $\Fp$-homology manifolds 
of the same dimension, let $G$ be 
a finite $p$-group and assume the following: The 
$\Fp$-cohomology of $M$ and $N$ is concentrated in even 
degrees and 
\[
    p > {\rm max} \{\dim_{\Fp} H^*(M) , \dim_{\Fp} H^*(N)\} \, .
\]
Let $f:M \ra N$ be $G$-equivariant. If the degree 
of $f$ is a nonzero number in $\Fp$, then the induced
map $M^G \ra N^G$ is surjective. 
\end{thm}

Note that by imposing these restrictions, for any 
normal subgroup $H < G$, the $\Fp$-cohomology 
of $M^H$ and $N^H$ will be concentrated in 
even degrees and all the spectral sequences
calculating equivariant cohomologies will 
collapse at the $E_1$-level.  Now  an inductive
proof can be carried out as indicated above.  

Without these strong additional assumption, we 
can prove the lifting assumption for differentiable actions
using the following general fact (cf.~\cite{Br}, Proposition
(2.3))  

\begin{lem} \label{abelian}
Let $G$ be a finite abelian $p$-group  and 
let $H < G$ be a subgroup. Let 
$E \ra X$ be a real linear $G$-bundle, if $p$ is odd,  or
a complex linear $G$-bundle, if $p=2$.  Each point 
in $X$ is assumed to have isotropy $H$ and $H$ is 
assumed to act freely on each fibre. Then the (nonequivariant) 
Euler class of $E$ lies in the image of 
\[
   H^*(X / (G/H)) \ra H^*(X) \, .
\]
\end{lem}

\begin{proof} Note that the claim is obvious, if $H$ is a direct
factor of $G$, because then the whole bundle carries a free 
$G/H$-action. In the following, we again concentrate on the 
special case $G = \Z / p^2$, $H= \Z /p$,  the general case being 
left to the reader. Let $\phi : G \times E \ra E$ be the given $G$-action on $E$. As an $H$-bundle, 
$E$ splits as a direct sum
\[
     \bigoplus_{j=1}^{p-1} E_j \otimes V_j
\]
where each $E_j$ is a real vector bundle and $V_1, \ldots , V_{p-1}$ 
are the disctinct nontrivial irreducible real $H$-representations. 
If $p$ is odd, each $V_j$ is a $1$-dimensional complex vector space
with $H$ acting complex linearly and the same is true in the 
case $p=2$ by our additional assumption.  It follows that 
the fibrewise $H$-action on $E$ can be propagated to 
a fibrewise $G$-action $\psi$. In particular, the $H$-action 
induced by $\psi$ coincides with the former fibrewise $H$-action
on $E$.

Now define a new $G$-action on $E$ by
\[
      (g, y ) \mapsto    \psi(g^{-1} ,   \phi(g , y)) \, .
\]
This action and the old $G$-action on $X$ define a new 
$G$-bundle $E' \ra X$. This bundle is nonequivariantly 
isomorphic to $E \ra X$, but carries an induced free $(G/H)$-action
which coincides with the old $(G/H)$-action on the base. 
\end{proof}

\begin{cor} Let $M$ and $N$ be oriented closed smooth 
$G$-manifolds, where $G= \Z /p^2$. If $p=2$, we 
assume additionally, that the normal bundle of $M^H$ in $M$ and 
of $N^H$ in $N$  are complex $G$-bundles. Then the above 
lifting assumption 
for  the $H$-equivariant Euler classes is fullfilled. 
\end{cor}

\begin{proof} The fixed point sets $M^H \subset M$ 
and $N^H \subset N$ have normal $G$-bundles $\nu^M \ra M^H$
and $\nu^N \ra N^H$. The equivariant Euler classes $e_H^M$ and 
$e_H^N$ are nothing but the usual Euler classes of 
the induced vector bundles 
\begin{eqnarray*}
   \nu^M \times_H EH & \ra &  M^H \times BH\, , \\ 
   \nu^N \times_H EH & \ra &  M^H \times BH
\end{eqnarray*}
Because $G$ is abelian, these bundles are again $G$-bundles (the 
$G$-operations are induced by the  $G$-operations on $\nu^M$ and $\nu^N$).  
After taking 
the cartesian product with the identity $EK \ra EK$, we get induced $G$-bundles
\begin{eqnarray*}
     EK \times (\nu^M \times_H EH)  & \ra & EK \times (M^H \times BH) \, , \\
     EK \times (\nu^N \times_H EH)  & \ra & EK \times (N^H \times BH)  
\end{eqnarray*}
where now $G$ is acting diagonally on the factor $EK$  (using the canonical map $G \ra K$) and the other factor. 
These $G$-vector bundles have the following properties \begin{itemize}
    \item Each point in the base has isotropy $H \subset G$.
    \item The group $H$ is acting freely on each fibre. 
\end{itemize}
Our claim now follows from Lemma \ref{abelian}.
\end{proof}

Based on this observation, it is not difficult to set up an inductive 
proof whose first steps were carried out above, in order 
to show once again the differentiable case of  Browder's theorem.

One could start trying to generalize Lemma \ref{abelian} in order 
to prove  generalizations of Browder's theorem similar
to those discussed in the first part of this paper. 
However, it seems hard to prove Lemma \ref{abelian} without the assumption 
of smoothness of the actions involved.

\end{document}